%% file: main2.tex
\documentclass[journal]{IEEEtran}

\ifCLASSINFOpdf
\fi

\usepackage{amsmath,amssymb,amsthm,amsfonts,mathrsfs}
\usepackage{rotating}

\newtheorem{theorem1}{Theorem}

\newtheorem{definition1}{Definition}

\newtheorem{lemma}{Lemma}
\newtheorem{assumption}{Assumption}
\newtheorem{corollary1}{Corollary}
\newtheorem{example}{Example}
\newtheorem{remark}{Remark}
\usepackage{enumerate}
\usepackage{verbatim}
\usepackage{pgf}
\usepackage{multicol}
\usepackage{booktabs}
\usepackage{multirow}
\usepackage{url}

\newcommand{\diag}{\mathop{\mathrm{diag}}}
\newcommand{\hvec}{\text{vec}}
\newcommand{\ngen}{n_G}
\newcommand{\real}{\text{Re}}
\newcommand{\imag}{\text{Im}}

\usepackage{ifthen}
\newboolean{showcomments}
\setboolean{showcomments}{true}
\newcommand{\eyb}[1]{\ifthenelse{\boolean{showcomments}}
{{ \color{red} {(Eilyan says:  #1)} }} {}}
\newcommand{\rl}[1]{  \ifthenelse{\boolean{showcomments}}
{{ \color{blue}{(Raphael says:  #1)} }}  {} }

\input{header.tex}

\usepackage[mathscr]{euscript}

\newcommand{\Ascr}{\mathcal{A}}

\newcommand{\Ab}{\mathbf{A}}
\newcommand{\Cb}{\mathbf{C}}
\newcommand{\Xb}{\mathbf{X}}
\newcommand{\Qb}{\mathbf{Q}}
\newcommand{\Lambdab}{\mathbf{\Lambda}}

\begin{document}

\title{Nondegeneracy and  Inexactness of Semidefinite Relaxations of Optimal Power Flow}

\author{Raphael~Louca$^1$,~\IEEEmembership{}
        Peter~Seiler$^2$,~\IEEEmembership{}
        and~Eilyan~Bitar$^1$.~\IEEEmembership{}

\thanks{1. R. Louca and E. Bitar are with the School of Electrical and Computer Engineering, Cornell University, Ithaca, NY 14853, USA.
        e-mail: {\tt\small rl553, eyb5@cornell.edu}}%
\thanks{2. P. Seiler is with the Department of Aerospace Engineering and Mechanics, University of Minnesota, Minneapolis, MN 55455, USA. e-mail:      {\tt\small seile017@umn.edu}}%
\thanks{Supported in part by NSF grant ECCS-1351621, NSF grant CNS-1239178, PSERC under sub-award S-52, US DOE under the CERTS initiative, and the University of Minnesota Institute on
the Environment under IREE Grant RL-0011-13.}}

\maketitle

\begin{abstract}
The Optimal Power Flow (OPF)  problem can be reformulated as a nonconvex Quadratically Constrained Quadratic Program (QCQP). There is a growing body of work on the use of semidefinite programming relaxations to solve OPF. The relaxation is exact if and only if the corresponding optimal solution set contains a rank-one matrix. In this  paper, we establish sufficient conditions guaranteeing  the nonexistence of a rank-one matrix in said optimal solution set. In particular, we show that under mild assumptions on problem nondegeneracy, any optimal solution to  the semidefinite relaxation will have rank greater than one, if the number of equality and active inequality constraints is at least twice the number of buses in the network. The sufficient condition holds for arbitrary network topologies (including tree networks). We empirically evaluate the practical implications of these results on several test cases from the literature.
\end{abstract}

\begin{IEEEkeywords}
Optimization, Semidefinite Programming, Optimal Power Flow, Nondegeneracy, Inexactness.
\end{IEEEkeywords}

\input{introduction.tex}

\input{problem_form.tex}

\input{main_results.tex}

\input{numerical_studies.tex}

\input{conclusion.tex}

\bibliographystyle{ieeetr}	
\bibliography{references_bib}

\input{appendix}

\end{document}

%% file: header.tex
\usepackage[margin=0pt,font=small,labelfont=bf,justification=justified,format=hang]{caption}

\usepackage{epsfig,amsfonts,latexsym,amssymb,shadow,amsmath,wrapfig}
\newcommand{\beq}{\begin{equation}}
\newcommand{\eeq}{\end{equation}}
\newcommand{\beqa}{\begin{eqnarray}}
\newcommand{\eeqa}{\end{eqnarray}}
\newcommand{\beqan}{\begin{eqnarray*}}
\newcommand{\eeqan}{\end{eqnarray*}}

\newcommand{\Rn}{{\mathbb R}^n}

\newcommand{\R}{{\mathbb R}}
\newcommand{\C}{{\mathbb C}}

\newcommand{\Cn}{\mathbb{C}^n}

\newcommand{\rank}{\mbox{rank}}

\newcommand{\bmat}[1]{ \begin{bmatrix} #1 \end{bmatrix}}

\newcommand{\Ical}{{\cal I}}
\newcommand{\Acal}{{\cal A}}
\newcommand{\Bcal}{{\cal B}}

\newcommand{\Dcal}{{\cal D}}
\newcommand{\Ecal}{{\cal E}}

\newcommand{\Hcal}{{\cal H}}

\newcommand{\Mcal}{{\cal M}}
\newcommand{\Ncal}{{\cal N}}
\newcommand{\Pcal}{{\cal P}}

\newcommand{\Rcal}{{\cal R}}

\newcommand{\Scal}{{\cal S}}

\newcommand{\white}[1]{{\color{white} #1}}

\newcounter{l1}
\newcounter{l2}
\newcounter{l3}
\newcommand{\bdotlist}{\begin{list}{$\bullet$}{}}
\newcommand{\bboxlist}{\begin{list}{$\Box$}{}}
\newcommand{\bbboxlist}{\begin{list}{\raisebox{.005in}{{\tiny
$\blacksquare$ \ \ }}}{}}
\newcommand{\bdashlist}{\begin{list}{$-$}{} }
\newcommand{\blist}{\begin{list}{}{} }
\newcommand{\barablist}{\begin{list}{\arabic{l1}}{\usecounter{l1}}}
\newcommand{\balphlist}{\begin{list}{(\alph{l2})}{\usecounter{l2}}}
\newcommand{\bAlphlist}{\begin{list}{\Alph{l2}.}{\usecounter{l2}}}
\newcommand{\bdiamlist}{\begin{list}{$\diamond$}{}}
\newcommand{\bromalist}{\begin{list}{(\roman{l3})}{\usecounter{l3}}}

\newcommand{\prf}[1]{ \noindent {\em Proof:} \, #1 \hfill $\blacksquare$}



%% file: introduction.tex
\section{Introduction}
\IEEEPARstart{T}{he} AC Optimal Power Flow Problem (OPF) is a nonlinear optimization problem where the objective is to minimize the total cost of generation subject to power balance constraints described by Kirchhoff's current and voltage laws
and operational constraints  reflecting bounds on real and reactive power generation, branch flows, and bus voltage magnitudes. OPF is a nonconvex optimization problem that is NP-hard in general. The problem nonconvexity derives from the nonconvex quadratic dependency of the feasible set on the set of complex bus voltages.
 
A recent stream of work has explored the application of semidefinite relaxations to solve OPF. Essentially, this convex relaxation entails the reformulation of OPF as a nonconvex Quadratically Constrained Quadratic Program (QCQP), which is then relaxed to a linear program over the cone of Hermitian positive semidefinite matrices. A detailed exposition on this relaxation approach can be found in \cite{low2014convex,low2014convex2} and the references therein.  The relaxation is said to be exact if and only if its
optimal solution set contains a rank-one matrix. In particular, several recent papers \cite{zhang2013geometry, bose2012quadratically, lavaei2012geometry,gan2013exact} have  established sufficient conditions on the OPF feasible region  under which the semidefinite relaxation is exact for networks with tree topology. One such sufficient condition amounts to the allowance of \emph{load over-satisfaction}. We refer the reader to Low  \cite{low2014convex2} for a comprehensive survey of these results.

While the importance of such exactness results over radial networks is self-evident, there remains an incomplete understanding of the behavior of semidefinite relaxations for more general network structures and OPF feasible regions. Of interest is the recent empirical observation that certain three-bus OPF problems yield inexact semidefinite relaxations \cite{lesieutre2011examining}. With this in mind, we pursue a deeper understanding of the limitations of the semidefinite relaxation technique for OPF. 

\vspace{-0.4cm}
\subsection{Contribution}

In what follows, we first exploit dual nondegeneracy of semidefinite programs -- as defined
by Alizadeh et al. in \cite{alizadeh1997complementarity} -- to establish \textit{uniqueness}  of
primal optimal solutions to complex semidefinite programs. 
Then, we establish \emph{necessary conditions}  for the
exactness of semidefinite relaxations of OPF for general
network structures. Or, interpreted in the contrapositive, we
offer sufficient conditions for the nonexistence of rank-one
optimal solutions to such relaxations. In particular, we show
that, under a mild assumption of
primal nondegeneracy, any optimal solution to the semidefinite
relaxation will have \emph{rank greater than one}, if the number of
equality and active inequality constraints is at least twice the
number of buses in the network. Of import is the fact that
this sufficient condition holds for arbitrary network topologies
(including trees). One would naturally expect this condition
to hold when the number of load buses in the network
is sufficiently large – as is common for traditional (radial)
distribution networks. And, indeed, our numerical analysis of
several test cases from the literature verifies this intuition.

\vspace{-0.4cm}
\subsection{Organization}
The remainder of the paper is organized as follows. 
In Section \ref{sec:prob_form}, we present the semidefinite relaxation of OPF  and describe assumptions and preliminary results required to establish our main result. Section \ref{sec:main_results} contains our main results. Theorem \ref{theorem:non_existence_opt} and Corollary \ref{coro:non_existence_opf} offer  sufficient conditions for the nonexistence of rank-one optimal solutions to the semidefinite relaxation of OPF.  In Section \ref{sec:numerical}, we perform a practical demonstration of our theory on several representative power system networks. Conclusions and directions for future research are given in Section \ref{sec:conclusion}.

%% file: problem_form.tex
\section{Problem Formulation} \label{sec:prob_form} 
\subsection{Notation}
Let $\R$ be the field of real numbers and $\C$ the field of complex numbers. For $z \in \C$, let $\real(z)$ and $\imag(z)$ be the real and imaginary parts of $z$ respectively.  For a matrix $A \in \C^{m \times n}$, let $a_{ij}$ be its $(i,j)$ entry. We denote by $A^\top$ and $A^*$ the transpose and complex conjugate transpose of $A$, respectively. Let $\Rn$ be the $n$-dimensional real vector space equipped with the usual inner product $\langle x, y \rangle := x^\top y$. We denote by $e_i$ the $i^{\text{th}}$ standard basis vector in $\Rn$. For square matrices $Y_1,\dots,Y_r$, we denote by $\diag(Y_1,\dots,Y_r)$, the square matrix with $Y_1,\dots,Y_r$ in the diagonal blocks and zeros in all other blocks. Let $\Hcal^n$ be the space of Hermitian $n \times n$ matrices, a real vector space of dimension $n^2$. 
We equip this space with the the Frobenius inner product, i.e., for $A, B \in \Hcal^n,$  $\smash{A \bullet B :=\sum_{i,j}a_{ij}b_{ij}^*}$. Let $\Hcal_+^n \ (\Hcal_{++}^n)$ be the space of Hermitian positive semidefinite (definite) matrices of order $n$. We write $A \succeq (\succ) B$, if and only if $A,B \in \Hcal^n$ and $A - B \in \Hcal^n_+ (\Hcal_{++}^n)$. 
We let  $|\Scal|$ be the cardinality of a set $\Scal$ and define set-theoretic difference as $\Scal_1 \setminus \Scal_2:=\{x \ | \ x \in \Scal_1 \text{ and } x \not\in \Scal_2\}.$

\subsection{Semidefinite Programming Relaxation}
Central to our analysis is the semidefinite relaxation of the OPF problem. Its derivation entails the exact reformulation
of the OPF problem as a semidefinite program subject to a rank-one inequality constraint on the set of feasible matrices \cite{low2014convex}. 
The semidefinite programming relaxation is obtained by removing the rank constraint and it can be written compactly as 
\begin{alignat} {6} \label{sdp_primal1}
& \underset{X \in \Hcal^n}{\text{minimize}}  
	 & & C \bullet X    \nonumber  \\
	 & \text{subject to}  \ \ 
	&& A_k \bullet X= b_k, \ \ \   \text{for all } k = 1,\dots,m,   \\
	&&& A_k \bullet X \le b_k, \ \ \  \text{for all } k = m+1,\dots,m+\ell,  \nonumber \\ 
	&&&  X \succeq 0, \nonumber
 \end{alignat}  
where $n$ is the number of nodes in the network, $C, A_1,\dots,$ $A_{m+\ell} \in \Hcal^n$ and $b:=[b_1,\dots,b_{m+\ell}]^\top \in \R^{m+\ell}$. The dual problem of (\ref{sdp_primal1}) is
\begin{equation}  \label{sdp_dual}
\begin{aligned}  
& \underset{y \in \R^{m+\ell}}{\text{maximize}}  
	& & -b^\top y \\ 
	&\text{subject to}  
	&& C + \sum_{k=1}^{m+\ell} y_k A_k \succeq 0, \\
	&&& y_k \ge 0, \qquad \qquad \quad    k = m+1,\dots,m+\ell.   
 \end{aligned}  
\end{equation}
Moving forward, we let $Z(y):= C + \sum_{k=1}^{m+\ell} y_kA_k$
and define the index sets $\Ecal:=\{1,\dots,m\} \text{ and } \Ical:=\{m+1,\dots,m+\ell\}$ for concision in notation. 
Let $\Pcal$ ($\Dcal$) be the feasible set of the primal (dual) problem and $\Pcal^\circ$ $(\Dcal^\circ)$ the corresponding optimal solution set of the primal (dual) problem. For any matrix $X \in \Pcal,$ let $$\Ascr(X):=\{k\in \Ical \ | \ A_k \bullet X =b_k\}$$ be the set of inequality constraints that are active at $X$ and define $a(X):=|\Ascr(X)|.$

\subsection{Assumptions}

We make the following assumptions about the primal-dual pair (\ref{sdp_primal1}) -- (\ref{sdp_dual}), which apply throughout the paper.

\begin{assumption}  \label{ass:ass1} \normalfont There exists an $X \in \Pcal$ such that $X \succ 0$ and $A_k \bullet X < b_k$ for all $k \in \Ical$ and a $y \in \Dcal$ such that $Z(y) \succ 0$ and $y_k > 0$ for all $k \in \Ical$. 
\end{assumption}
Assumption  \ref{ass:ass1} is a Slater Condition. It guarantees strong duality to hold. This in turn implies
 a complementarity condition \cite{alizadeh1997complementarity}. Namely, for any pair of primal-dual optimal solutions $(X,y) \in \Pcal^\circ \times \Dcal^\circ$, it holds that $\rank(X)+\rank(Z(y)) \le n$.

\begin{assumption} \normalfont \label{ass:ass2} Strict complementarity holds between any pair of primal-dual optimal solutions.
\end{assumption}
A primal-dual optimal solution pair $(X,y) \in \Pcal^\circ \times \Dcal^\circ$ is said to satisfy \textit{strict complemnetarity} if $\rank(X)+\rank(Z(y))=n$. The importance of strict complementarity will be made apparent in the sequel (c.f. Theorem \ref{theorem:non_existence_opt}).

Several authors  \cite{alizadeh1997complementarity,pataki2001generic} have shown that strict complementarity is a \textit{generic property} of optimal solutions to semidefinite programs. A property of semidefinite programs is \textit{generic} \cite{alizadeh1997complementarity} if it holds for almost all problem instances $(C,b,A_1,\dots,A_{m+\ell})$; that is, the set of problem instances for which said property fails to hold has measure zero.

\begin{assumption} \normalfont \label{ass:ass3} The  matrices in the set $\{A_k \ | \  k \in \Ecal\}$ are linearly independent.
\end{assumption} 
The linear independence assumption is without loss of generality.  If the matrices  $A_k$, $k=1,\dots,m$ are linearly dependent, one can choose a basis of, say $p<m$, matrices from $\{A_k \ | \  k \in \Ecal\}$ and remove the other $m-p$ equality constraints to establish an equivalent problem in which the assumption holds.
The significance of this assumption will be made apparent in the sequel when we define primal and dual nondegeneracy in  Definitions \ref{def:nondegeneracy_primal}--\ref{def:nondegeneracy_dual}.
In Appendix \ref{app:Linear_Inde}, we introduce a linear transformation of the problem data to check for linear independence of a set of Hermitian matrices. 

In the following Section \ref{sec:IIC}, we offer a simple extension of Alizadeh's linear algebraic characterization of \textit{nondegeneracy} of solutions to semidefinite programs defined over the real symmetric positive semidefinite cone \cite{alizadeh1997complementarity,alizadeh1997optimization}
to the Hermitian positive semidefinite cone.

\vspace{-0.3cm}
\subsection{Definitions and Related Results} \label{sec:IIC}
Of central importance to our results is the notion of \textit{nondegeneracy} of solutions to semidefinite programs. In \cite{alizadeh1997complementarity}, the authors define nondegeneracy of primal and dual feasible points for semidefinite programs in standard form. In \cite{nayakkankuppam1999conditioning}, the authors extend this analysis to block-structured semidefinite programs, where the decision matrix is restricted to have block diagonal structure. In Appendix \ref{sec:AppII}, we show how to reformulate problem (\ref{sdp_primal1}) as a block semidefinite program. We thus have the following definition, which is a special case of the definition of primal nondegeneracy for block-structured semidefinite programs in \cite{nayakkankuppam1999conditioning}.

\begin{definition1}[Primal Nondegeneracy] \label{def:nondegeneracy_primal} \normalfont 
Let $X \in \Pcal$ and suppose that $\rank(X) = r$. Let $X = Q \Lambda Q^*$ be an eigenvalue decomposition of $X$, where $\Lambda = \diag(\lambda_1,\dots,\lambda_r,0,\dots,0) \in \R^{n \times n}$ and $Q \in \C^{n \times n}$. Partition $Q$ as $Q = \smash{\bmat{Q_1 & Q_2}},$ where $Q_1 \in \C^{n \times r}$ and $Q_2 \in \C^{n \times (n-r)}$, and define the matrices
$$	B_k^p := \bmat{Q_1^*A_kQ_1 & Q_1^* A_k Q_2 \\ Q_2^*A_kQ_1 & 0} \ \  \text{for all} \  \ k\in \Ecal \cup \Ascr(X).$$
Then $X $ is \emph{primal nondegenerate}, if and only if the matrices $\{B_k^p \ | \ k\in \Ecal \cup \Ascr(X)\}$ are linearly independent in $\Hcal^n$. \qed
\end{definition1}

We similarly define nondegeneracy of dual feasible points.

\begin{definition1}[Dual Nondegeneracy] \label{def:nondegeneracy_dual} \normalfont
Let $y \in \Dcal$ and suppose that $\rank(Z(y))=s.$ Let $Z(y) = P \Sigma P^*$ be an eigenvalue decomposition of $Z(y)$, where $\Sigma = \diag(0,\dots,$ $0,\sigma_{n-s+1},\dots,\sigma_n) \in \R^{n \times n}$ and $P \in \C^{n \times n}$. Partition $P$ as $P = \smash{\bmat{P_1 & P_2}},$ where $P_1 \in \C^{n \times (n-s)}$ and  $P_2 \in \C^{n \times s}$,  and define the matrices
$$ B_k^d := \bmat{P_1^*A_kP_1} \ \  \text{for all} \ \  k\in \Ecal \cup \Ascr(X). $$	
Then $y$ is \emph{dual nondegenerate}, if and only if  the  matrices $\{B_k^d \ | \ k\in \Ecal \cup \Ascr(X)\}$ span $\Hcal^{n-s}$. \qed
\end{definition1}

Note that if $(X,y)\in \Pcal^\circ \times \Dcal^\circ$ is a primal-dual optimal solution pair that satisfies strict complementarity, then $Q_1=P_1$ and $Q_2=P_2$. 

\begin{remark}[Transversality] \normalfont
Primal nondegeneracy has the following geometric interpretation. Let $\Mcal_r:=\{X\in \Hcal^n \ | \ \rank(X)=r\}$ be the set of Hermitian matrices of order $n$ that have rank $r$. A primal feasible point $X \in \Pcal \cap \Mcal_r$ is nondegenerate if and only if the orthogonal complement of the subspace spanned by the matrices $A_k, \ k \in \Ecal \cup \Ascr(X),$  intersects the tangent space to $\Mcal_r$ at $X$ transversally. See \cite{alizadeh1997complementarity, shapiro1997first} for  a definition of transversality. \qed
\end{remark}

Primal and dual nondegeneracy is related to uniqueness of optimal solutions to semidefinite programs. In particular, we have the following result from Alizadeh et al. \cite{alizadeh1997complementarity}. 

\begin{theorem1}[Uniqueness of Optimal Solutions] \normalfont \label{thm:dual_nond_unique_primal_optimal} 
Let $y \in \Dcal^\circ$ be dual nondegenerate. Then, there exists a unique primal optimal solution. Similarly, if $X \in \Pcal^\circ$ is primal nondegenerate, then there exists a unique dual optimal solution. \qed
\end{theorem1} 

\begin{remark} \normalfont
The assumption of dual nondegeneracy is natural, as it holds \emph{generically} for semidefinite programs \cite{alizadeh1997complementarity}. Under this assumption, if $X \in \Pcal^\circ$ has $\rank(X)>1$, then said  semidefinite relaxation of OPF is necessarily \emph{inexact}.   \qed
\end{remark}

%% file: main_results.tex
\section{Main Results} \label{sec:main_results}
In this section we leverage on primal nondegeneracy to establish sufficient conditions guaranteeing the inexactness  of  semidefinite relaxations of OPF.

\subsection{Nonexistence of Rank-One Solutions} \label{sec:main_results_nonexistence}

Beyond yielding uniqueness of dual optimal solutions to semidefinite programs, primal nondegeneracy can be exploited to provide a sufficient condition for the nonexistence of rank-one optimal solutions to semidefinite relaxations of OPF. More precisely, we have the following result, which follows readily from Definition \ref{def:nondegeneracy_primal}.

\begin{lemma} \normalfont \label{lem:aposteriori}
Let $X \in \Pcal$ be a primal nondegenerate feasible point. If 
$m + a(X) \ge 2n,$
then $\rank(X) > 1.$ \qed
\end{lemma}

\prf{
 Let $X \in \Pcal$ be primal nondegenerate and suppose that $\rank(X) = r$. Primal nondegeneracy of $X$ implies that the collection of matrices $B_k^p, \ k\in \Ecal \cup \Ascr(X)$ specified in Definition \ref{def:nondegeneracy_primal} are linearly independent in $\Hcal^n$. Hence, they span an $m+a(X)$ dimensional subspace in $\Hcal^n$. Recall that $\Hcal^n$  forms a vector space over the real numbers of dimension $n^2$. And let ${\hvec:\Hcal^n \rightarrow \R^{n^2}}$ be a function that maps Hermitian matrices of order $n$ to real vectors of length $n^2$ as follows:
\begin{equation*}
\begin{split}
	 \hvec(A) := [&a_{11}, \real(a_{21}),  \imag(a_{21}), \dots,\real(a_{n1}), \imag(a_{n1}),  a_{22}, \\ & \real(a_{23}),\imag(a_{23}),\dots,a_{nn}]^\top. 
	 \end{split}
	 \end{equation*}
Upon application of this transformation to $B_k^p$, it is straightforward to see that $$\hvec(B_k^p) = \bmat{\xi_k \\ 0}^\top \in \R^{n^2},$$ for all $k \in \Ecal \cup \Ascr(X)$, where $\xi_k \in \R^{n^2-(n-r)^2}$ and $0 \in \R^{(n-r)^2}$ is the zero vector. And, since the matrices $B_k^p, \ k\in \Ecal \cup \Acal(X)$ span an $m+a(X)$ dimensional space in $\Hcal^n$, it must be true that $m+a(X) \le n^2 - (n-r)^2$. It follows that $m+a(X) \le n^2 - (n-r)^2$, if and only if $r\ge n-[n^2-m-a(X)]^{1/2}.$
And finally, we have $n-[n^2-m-a(X)]^{1/2}> 1$, if and only if $m+a(X) \ge 2n$, thus completing the proof.
}

\begin{theorem1} \label{theorem:non_existence_opt} \normalfont Let $X^\circ \in \Pcal^\circ$ be a primal nondegenerate optimal solution.  If $m + a(X^\circ) \ge 2n$, then $\rank(X) > 1$ for all $X \in \Pcal^\circ$.  \qed
\end{theorem1} \vspace{-0.5cm}
\begin{proof}
	Since $X^\circ \in \Pcal^\circ$ is primal nondegenerate and optimal, there exists a unique dual optimal solution by Theorem \ref{thm:dual_nond_unique_primal_optimal}. Let $y\in \Dcal^\circ$ be this solution. Since $m + a(X^\circ) \ge 2n,$ it follows by Lemma \ref{lem:aposteriori} that $\rank(X^\circ) >1$ and by complementarity we must have that $\rank(Z(y)) \le n - \rank(X^\circ) < n-1$. Since strict complementarity is assumed to hold for every primal-dual optimal pair (c.f. Assumption \ref{ass:ass2}) and $y \in \Dcal^\circ$ is unique, it must be true that $\rank(X)>1$ for all $X \in \Pcal^\circ.$
\end{proof}
Theorem \ref{theorem:non_existence_opt} offers a sufficient condition for the nonexistence of rank-one optimal solutions. This condition can be verified \emph{a priori} in the event that the number of native equality constraints $m$ is large. We have a trivial Corollary of Theorem \ref{theorem:non_existence_opt} in the following result.

\begin{corollary1} \label{coro:coro_main_thm}\normalfont 
Suppose there exists an $X^\circ \in \Pcal^\circ$  that is primal nondegenerate. If $m \ge 2n$, then $\rank(X)>1$  for all $X \in \Pcal^\circ$. \qed
\end{corollary1}  

\begin{remark}[Exactness] \normalfont 
Naturally, the sufficient condition in Theorem \ref{theorem:non_existence_opt} and Corollary \ref{coro:non_existence_opf}   can be reformulated in the contrapositive as \emph{neccessary conditions for the exactness of  semidefinite relaxations of OPF.}  \qed
\end{remark}

\subsection{Implications on Semidefinite Relaxations of OPF}
In Section \ref{sec:main_results_nonexistence}, we showed that under mild assumptions on problem nondegeneracy,   optimal solutions to semidefinite programs will have  rank greater than one, if the number of equality and active inequality constraints is greater than or equal to twice the ambient dimension of the problem. In what follows, we discuss implications of this result to semidefinite relaxations of OPF. Namely, we specify sufficient conditions under which the semidefinite relaxation of OPF is guaranteed to fail.  We make this precise in what follows.

Consider a power network consisting of $n$ buses and let $\Ncal :=\{1,\dots,n\}$. Let $\Ncal_G \subseteq \Ncal$ denote the set of buses connected to
 generators and define $\ngen:=|\Ncal_G|$. 
 An essential constraint required by OPF is network power balance. The power balance equations, described by Kirchhoff's current and voltage laws, govern the relationship between complex bus voltages and  power injections. 
 More precisely, let $V \in \C^{n}$ be the vector of bus voltage phasors and $Y \in \C^{n \times n}$  the network admittance matrix. 
 Denote by $S_G \in \Cn$ and $S_D \in \Cn$ the vectors of complex bus power generation and demand, respectively. Note that $S_{G_i} = 0 $  for all $i \notin \Ncal_G$. 
 
The power balance equations can thus be expressed as 
\begin{equation} \label{eq:power_balance}
\begin{aligned} 
  S_{D_i} + Y^*e_ie_i^{\top} \bullet VV^*  =  0 \quad 
\end{aligned}
\end{equation}
for all $i \in \Ncal \backslash \Ncal_G$
 non-generator (i.e., load) buses and
\begin{equation}
 \begin{aligned}
 &P_{G_i}^{\min} \le \real\{S_{D_i} + Y^*e_ie_i^\top \bullet V V^* \} \le P_{G_i}^{\max} \\
 &Q_{G_i}^{\min} \le \imag \{S_{D_i} + Y^*e_ie_i^\top \bullet V V^* \} \le Q_{G_i}^{\max}
  \end{aligned} 
\end{equation} 
for all $i \in \Ncal_G$ generator buses. Here, $P_{G_i}^{\max}$ and $ Q_{G_i}^{\max}$ denote upper bounds on real and reactive power generation respectively. Similarly, $P_{G_i}^{\min}$ and $Q_{G_i}^{\min}$ denote lower bounds on real and reactive power generation, respectively.
Clearly, the power balance equation \eqref{eq:power_balance} amounts to $2(n - n_G)$ equality constraints. This coincides with twice the number of  load buses in the network. 
In addition to requiring power balance, it is common for OPF problems  to impose equality constraints on the bus voltage magnitudes of the form
\begin{equation} \label{eq:volt_mag}
\begin{aligned} 
  e_ie_i^{\top} \bullet VV^* = \overline{V}_i^2 
\end{aligned}
\end{equation}
for $i \in \Ncal_V$, where $\Ncal_V$ denotes the set of buses requiring fixed voltage magnitudes and $n_V :=|\Ncal_V|$. Here, $\overline{V}_i$ denotes the required voltage magnitude at bus $i \in \Ncal_V$. It follows that the number of equality constraints, $m$, inherent to the semidefinite relaxation of OPF will be at least
\begin{align} \label{eq:lower_bound_eq}
m \ge 2(n-\ngen) + n_V.
\end{align}
This implies an immediate corollary of Theorem \ref{theorem:non_existence_opt}, which we state without proof. 

\begin{corollary1} \label{coro:non_existence_opf} \normalfont  Suppose that $X^\circ \in \Pcal^\circ$ is a primal nondegenerate solution to the semidefinite relaxation of OPF. If $$ \ngen  \leq \frac{1}{2} ( a(X^\circ)  + n_V),$$ then $\rank(X) > 1$ for all $X \in \Pcal^\circ$.  In particular, if $n_G \leq \frac{1}{2} n_V,$ then $\rank(X)>1$ for all $X \in \Pcal^\circ.$ \qed
\end{corollary1} 

Clearly, Corollary \ref{coro:non_existence_opf} implies that OPF problems with a sufficiently small number of generator buses will fail to yield to semidefinite relaxations having rank-one optimal solutions. 

\begin{remark} \normalfont
Often, the voltage magnitude constraints enter as inequalities specifying an acceptable range of values. In practice, however, the lower and upper bounds on bus voltage
magnitudes are chosen to be close to 1 per unit for all buses $i \in \Ncal$,  because of strict requirements on power quality. As such, it is not uncommon to observe binding voltage magnitude constraints at optimality. \qed
\end{remark}

\begin{figure*}[t] 
\begin{center}
\includegraphics[width=2.25in]{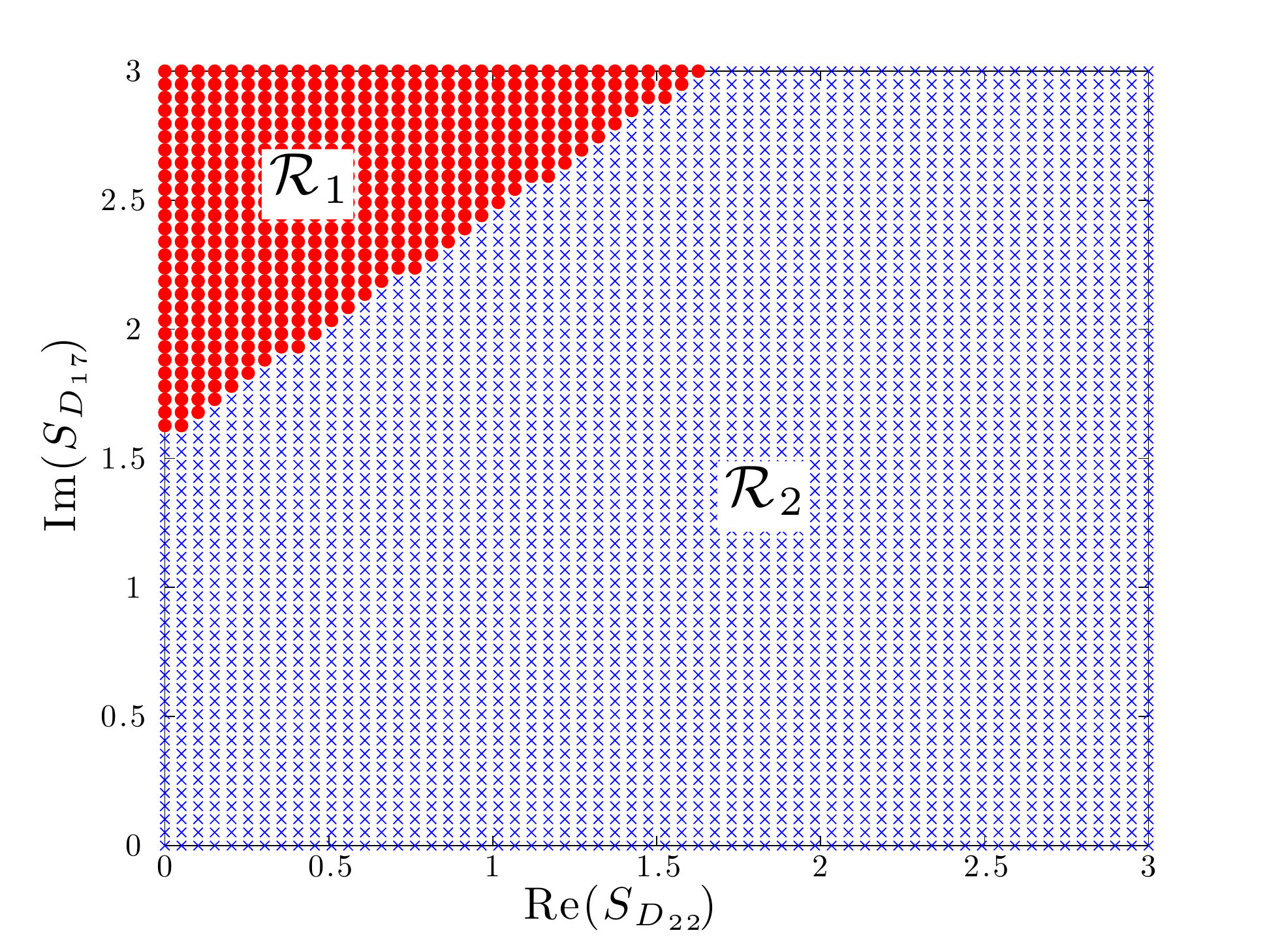}
\includegraphics[width=2.25in]{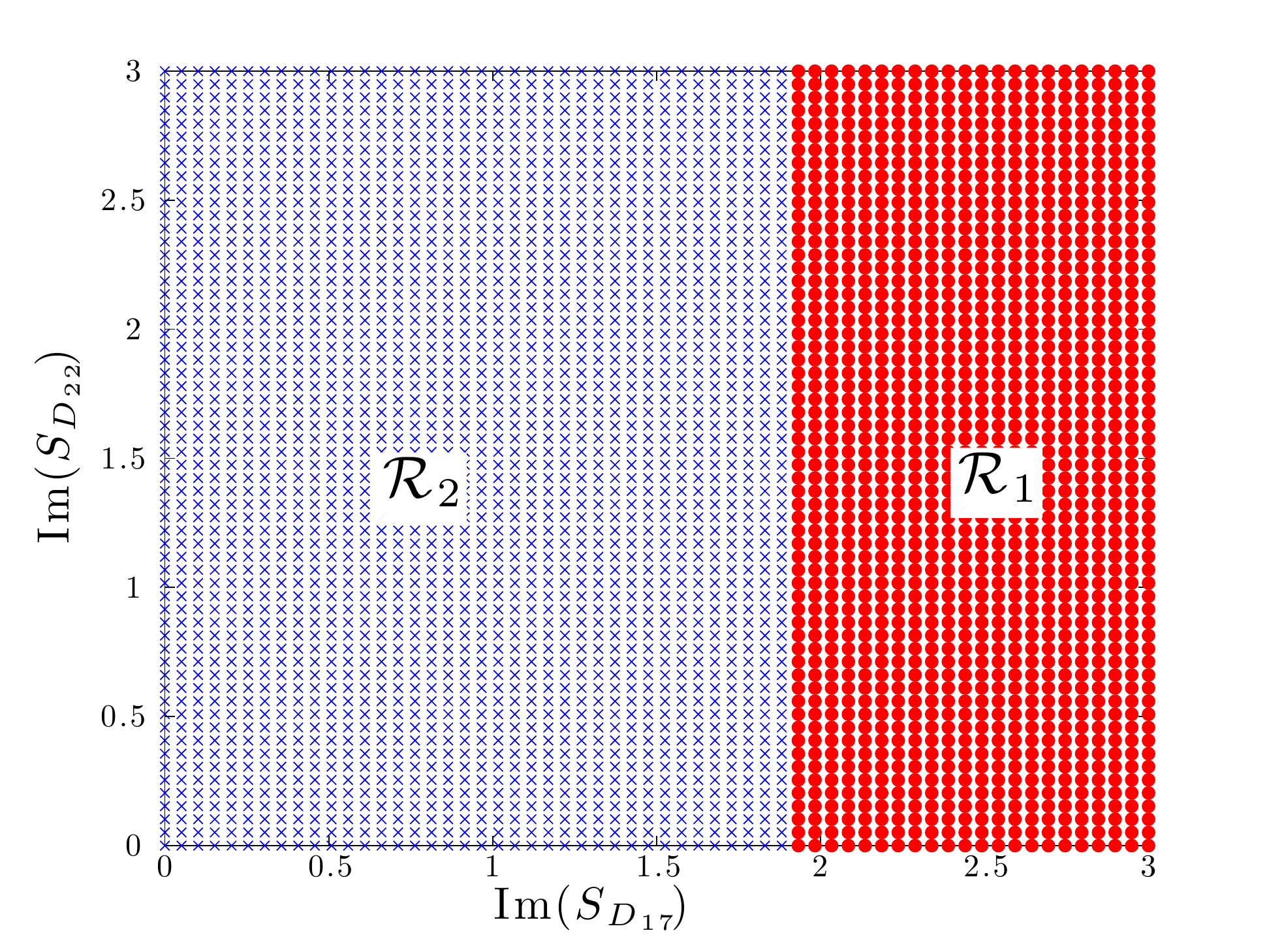}
\includegraphics[width=2.25in]{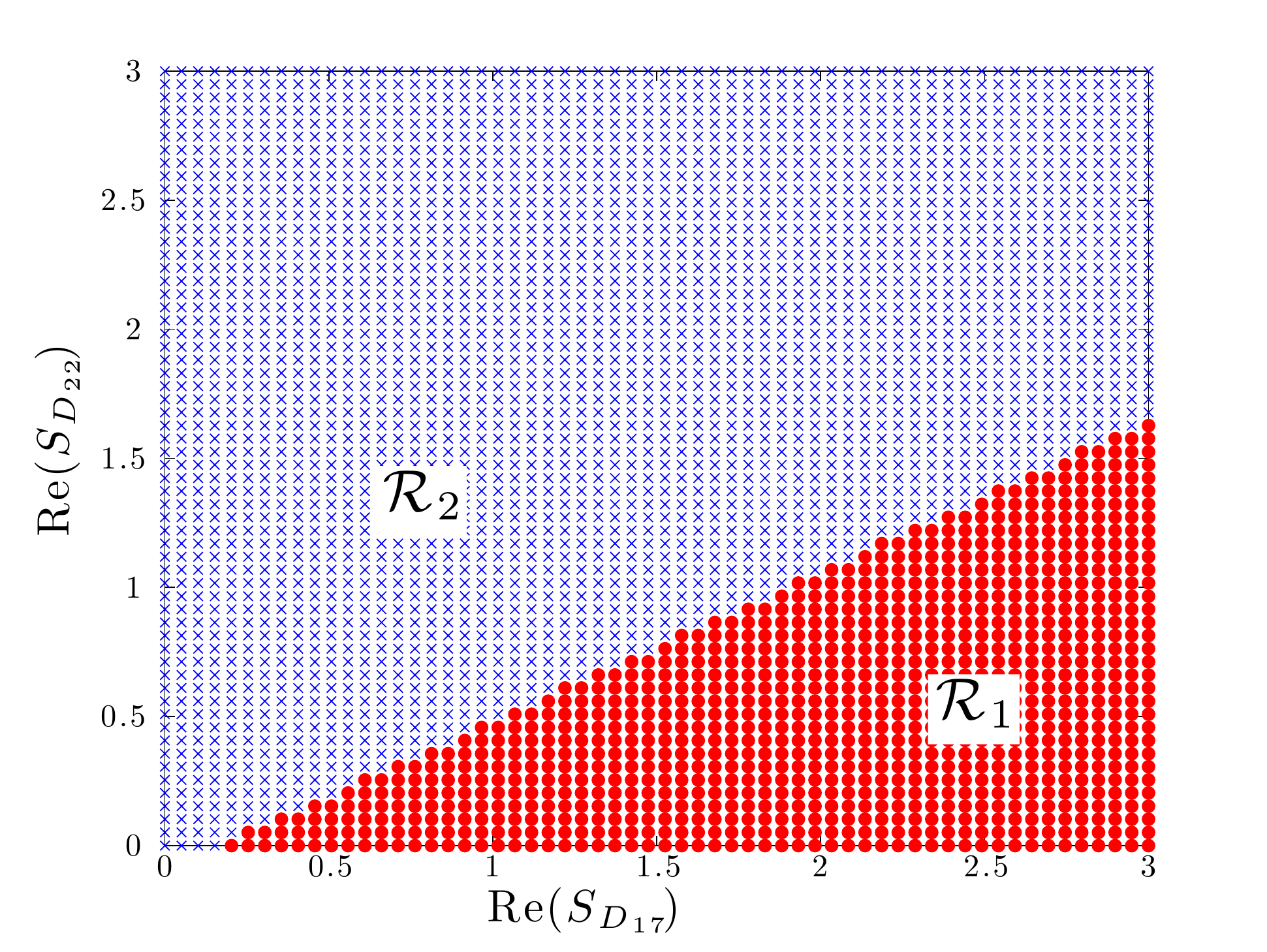}
\caption{IEEE 34-Node Test Feeder (Radial Network): Two dimensional sections of the sets of load profiles for which the semidefinite relaxation admitis unique rank-one ($\Rcal_1$) and unique high rank ($\Rcal_2$) optimal solutions. The red circles correspond to points in the set $\Rcal_1:=\{ S_D \in \C^{n} \ | \ (X,y)\in \Pcal^\circ \times \Dcal^\circ, \  \rank(X)=1, \ y \text{ is dual nondegenerate}\}$ and the blue crosses to points in the set $\Rcal_2:=\{ S_D \in \C^{n} \ | \ (X,y)\in \Pcal^\circ \times \Dcal^\circ, \  \rank(X)>1, \ y \text{ is dual nondegenerate}\}$. } 
\label{fig:empirical} 
\end{center}
\end{figure*}

\begin{remark}[Load-Oversatisfaction] \label{rem:load_over} \normalfont A number of sufficient conditions in the  literature guaranteeing  exactness of solutions to semidefinite relaxations of OPF over \textit{radial networks} rely on the so-called assumption of \textit{load-oversatisfaction} \cite{low2014convex2, sojoudi2011network}. 
In part, the assumption of load-oversatisfaction amounts to relaxing the power balance equations (at load buses) in (\ref{eq:power_balance}) to inequalities, where the complex power delivered to each node is allowed to exceed the power demanded. Under such assumption, it can be shown that semdefinite relaxations over radial networks will have optimal solutions of rank one. Moreover, it is claimed in \cite{sojoudi2011network} that such solutions will in general tend to satisfy the power balance equations, despite the allowance of load-oversatisfaction.
However, if the solution to the semidefinite relaxation has rank greater than one and the corresponding dual optimal solution is nondegenerate, then  the solution to the semidefinite relaxation when load oversatisfaction is allowed will \textit{necessarily violate} the power balance equations (c.f Theorem \ref{thm:dual_nond_unique_primal_optimal}). A similar argument can be made (for primal nondegenerate solutions),  when the number of generator buses is no greater than half the number of buses with fixed voltage magnitudes (c.f. Corollary \ref{coro:non_existence_opf}).\qed
\end{remark}

\begin{remark}[Optimal Voltage Regulation] \normalfont The optimal voltage regulation problem considered in \cite{lam2012optimal} amounts to an OPF problem over a radial network, where the bus voltage magnitudes are fixed at every bus in the network (i.e., $n_V = n$). Corollary \ref{coro:non_existence_opf} implies that any primal nondegenerate  solution to the corresponding semidefinite relaxation will have rank strictly greater than one, if $n_G \leq (1/2) n$. 
Equivalently, there exists a primal nondegenerate solution that has rank equal to one, only if $\ngen > (1/2)n$. Of critical import is the fact that traditional radial (distribution) networks possess few points of controllable generation (i.e., $\ngen   \ll n$).  
 This will likely change in the future, however, as the power system evolves to incorporate increased distributed generation (e.g. rooftop solar, plugin electrical vehicles, and controllable appliances).
\qed
\end{remark}

%% file: numerical_studies.tex
\section{Numerical Studies} \label{sec:numerical}

  We now empirically evaluate our theoretical results  on several test cases from the literature.
   Throughout this section, we consider linear objective functions of the form
$$ f_i(S_{G_i}) := c_{i_1}\real(S_{G_{i}}) + c_{i_0},$$
where $i \in \Ncal_G$ and $c_{i_1},c_{i_0}$ are non-negative scalars.

\subsection{Cyclic Networks with Unique High Rank Solutions}

We summarize, in Table \ref{table:case1}, several OPF test cases from the literature \cite{lesieutre2011examining, zimmerman2011matpower, scotland_cases} whose corresponding semidefinite relaxations are guaranteed  to be inexact according to Theorem \ref{thm:dual_nond_unique_primal_optimal}.
Namely, their semidefinite relaxations admit \emph{unique high rank optimal solutions}. The first four examples reveal the nonnecessity of the sufficient condition for inexactness outlined in Theorem \ref{theorem:non_existence_opt}, as $m + a(X) < 2n$ in each case. 
The  OPF case (14C)  in Table \ref{table:case1} is of
particular interest. This example\footnote{\label{footnotemark1} We provide a precise description of the data in \url{http://foie.ece.cornell.edu/~louca/opf/} } is a modified IEEE 14 Bus system for which the solution to the semidefinite relaxation is high rank and  satisfies the sufficient conditions of Theorem \ref{theorem:non_existence_opt} and Corollary \ref{coro:non_existence_opf}. Other examples in the literature that yield high rank optimal solutions include the 9 and 30 bus networks in \cite{zimmerman2011matpower}. For these OPF problems, however, the solution to the semidefinite relaxation is \textit{dual degenerate}. By adding a small resistance (e.g. $10^{-5}$ ohms) to a subset of the lines  with zero resistance, we obtain a dual nondegenerate optimal solution and a (unique) rank-one primal optimal solution 
 whose cost is within 0.002\% of the optimal value in the degenerate case. 

\white{\footnote{A small resistance of $10^{-5}$ was added to all lines having zero resistance.}}
  
\begin{table}[t]
\caption{Power system examples with semidefinite relaxations having \textit{unique} high rank solutions.} 
\label{table:case1} \centering%
\setlength{\tabcolsep}{5.5pt}
\begin {tabular*}{0.46\textwidth}{ccccccc}
\toprule
System &  \multirow{2}{*}{Ref.} & \multirow{2}{*}{$\rank(X)$} & $X \in \Pcal^\circ$  & $y \in \Dcal^\circ$  & \multirow{2}{*}{$m$ } & \multirow{2}{*}{$a(X)$ } \\ 
$(n)$ & &  & Nondege. & Nondege. &  &      \\ \toprule
3 & \cite{lesieutre2011examining}			    &  2  & $\surd$ &  $\surd$    	  & 1  & 4   		\\
5 & \cite{scotland_cases}  	    &  2  &  $\surd$ &    $\surd$     &  6  & 3   \\
39$^2$ & \cite{zimmerman2011matpower}		     & 2    & $\surd$ & $\surd$        & 58  & 14         \\ 
118  & \cite{zimmerman2011matpower}           &  2   & $\surd$ &  $\surd$     & 128  & 73       \\ 
14C$^1$ & --  &  2   & $\surd$    &   $\surd$    &  18   &   12      \\    \bottomrule
\end{tabular*}\vspace{-0.3cm}
\end{table}	  
\begin{figure}[h] \vspace{-0.8cm}
\centering
	\includegraphics[scale=0.57]{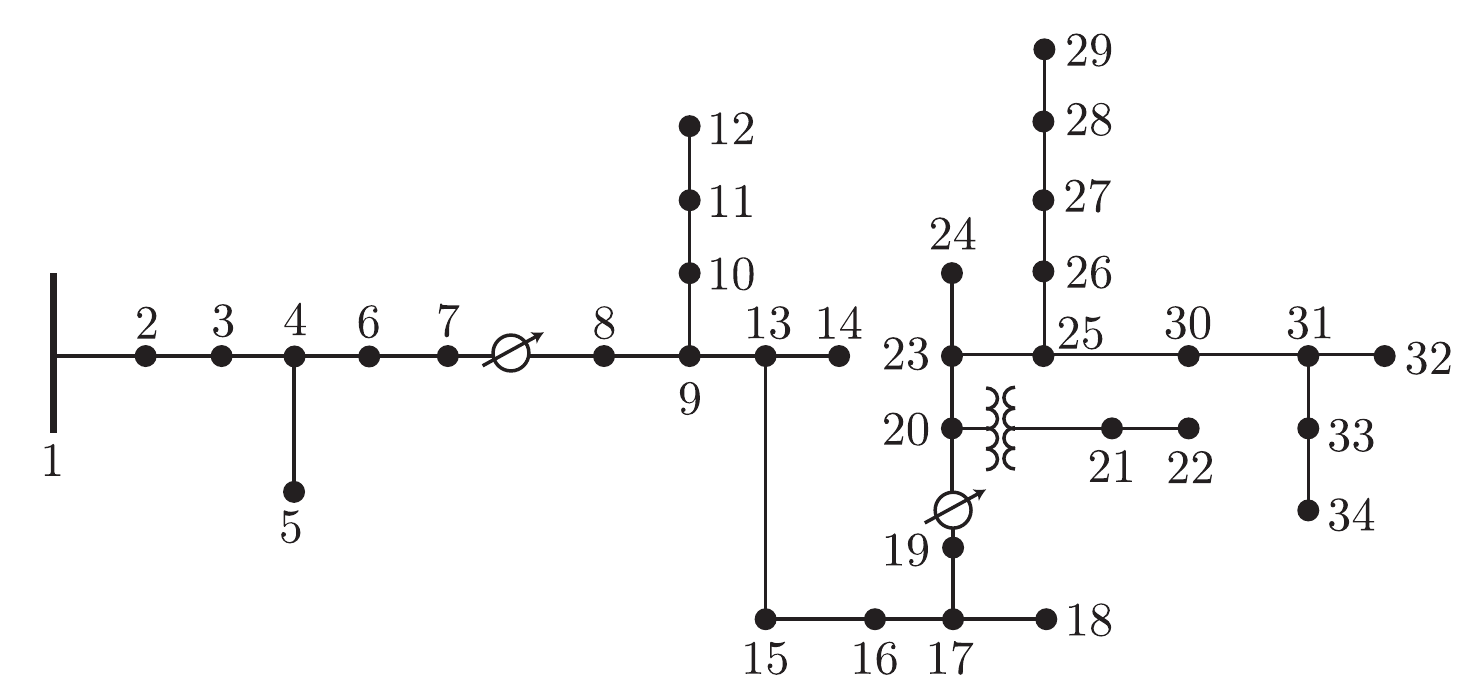} 
	\caption{IEEE 34 Node Test Feeder (Radial Network) }
	\label{fig:34bus} \vspace{-0.4cm}
\end{figure}

\vspace{-1.1cm}
\subsection{A Radial Network with Inexact Semidefinite Relaxation} 
We now empirically explore the extent to which semidefinite relaxations over radial networks prove to be inexact.
In particular, we consider the IEEE 34 bus test feeder case$^1$, which has a radial (tree) structure. Figure \ref{fig:34bus} offers a graphical illustration of the network topology. The OPF problem  considered consists of $n_G = 33$ generator buses and has $m = 35$ and $\ell = 134$ native equality and inequality constraints, respectively, describing its feasible set. \emph{Our objective is to evaluate the effect of the underlying load profile on the success of the corresponding semidefinite relaxation.}
We do so by conducting a parametric analysis over a subset of the system's load values. More precisely, we let $\real\{S_{D_i}\}$  and $\imag\{S_{D_i}\}$ range over the interval $[0, 3]$ \text{MW}  for buses $i \in \{17,22\}$.
As a practical consideration, we uniformly grid this parameter space $[0,3]^4$ by generating 10 equally spaced points in  $[0,3]$ along each coordinate axis.
We fix the remaining loads $S_{D_i}$,  at buses $i \in \Ncal \setminus \{17,22\}$, to their nominal values.
One can numerically verify (using Matpower \cite{zimmerman2011matpower}) that the OPF problem induced under each load profile in the parameter space considered has a nonempty feasible set.

We now proceed in describing our empirical observations. For each load profile considered, the resulting semidefinite relaxation yielded a dual nondegenerate optimal solution, which implies  uniqueness of the corresponding primal optimal solution. 
 Moreover, $99.6 \%$  of the load profiles considered produced primal nondegenerate optimal solutions. This is reassuring, as primal and dual nondegeneracy are \textit{generic} properties of semidefinite programs \cite{alizadeh1997complementarity}. Of interest is our empirical observation that $93.69\%$ of the load profiles considered yielded unique primal optimal solutions of high rank, while only $6.31\%$ admitted rank-one optimal solutions. As such, the semdefinite relaxation is \emph{inexact} for the majority of load profiles considered. 
 In addition,  all primal nondegenerate optimal solutions $X$ of high rank were observed to satisfy the sufficient condition for inexactness specified in Theorem \ref{theorem:non_existence_opt}, i.e., $m+a(X) \ge 2n$. 
 And all primal nondegenerate optimal solutions $X$ of rank one satisfied $m+a(X)=2n-1$.

In Figure \ref{fig:empirical}, we visualize two-dimensional sections of the sets of load profiles for which the semidefinite relaxation admits unique rank-one ($\Rcal_1$) and unique high rank ($\Rcal_2$) optimal solutions.
 It is important to note that the set $\Rcal_2$ specifies the set of loads for which a solution obtained from a semidefinite relaxation \textit{allowing load oversatisfaction} is guaranteed to violate the power balance equations (c.f. Remark \ref{rem:load_over}).

\begin{remark} \normalfont While our example highlights a radial network for which the semidefinite relaxation is inexact for a large fraction of load profiles considered, one should take care not to blindly apply such conclusions to the more general family of all radial networks. Rather, this example is meant to reveal that semidefinite relaxations over radial networks can indeed fail over practical operating regimes.
This should serve to 
stimulate additional research focused on the refinement of existing necessary and/or sufficient conditions for the 
exactness of semidefinite relaxations over networks with tree topological structure.\qed \end{remark}

\vspace{-0.7cm}

%% file: conclusion.tex
\section{Conclusion} \label{sec:conclusion}
This paper considered the nonconvex Optimal Power Flow (OPF) problem and the corresponding semidefinite relaxation. By leveraging on the theory of nondegeneracy in  semidefinite programming, we construct sufficient conditions under which solutions to semidefinite relaxations of OPF over arbitrary network topologies are guaranteed to have rank strictly greater than than one. This condition is shown to hold when  the number of  equality constraints inherent to the OPF problem is sufficiently large. 
\vspace{-0.4cm}

%% file: appendix.tex
\appendices
\vspace{-0.3cm}
\section{Linear Independence of Hermitian Matrices} \label{app:Linear_Inde}
For $A \in \Hcal^n$, let $\hvec:\Hcal^n \rightarrow \R^{n^2}$ be a function that maps Hermitian matrices of order $n$ to column vectors of length $n^2$ as follows:
\begin{small}
\begin{equation*}
\begin{split}
	 \hvec(A) := [&a_{11}, \real(a_{21}),  \imag(a_{21}), \dots,\real(a_{n1}), \imag(a_{n1}),  a_{22}, \\ & \real(a_{23}),\imag(a_{23}),\dots,a_{nn}]^\top. 
	 \end{split}
	 \end{equation*}
\end{small}{\raggedright
Let $V_1,\dots,V_d \in \Hcal^n$ be Hermitian} matrices of order $n$. Then, $V_1,\dots,V_d$ are linearly independent if and only if the matrix
$$V:=[\hvec(V_1),\dots, \hvec(V_d)] \in \R^{n^2 \times d}$$ has full rank.

\vspace{-0.5cm}
\section{Problem Reformulation} \label{sec:AppII}
Let $\Bcal^{n+\ell}$ denote the space of all Hermitian $(n+\ell) \times (n+\ell)$ block diagonal matrices, consisting of one $n$-dimensional diagonal block and $\ell$ one-dimensional diagonal blocks. We refer to the $i^{\text{th}}$ diagonal block of a matrix $\Xb \in \Bcal^{n+\ell}$ as $\Xb(i)$. The dimension of the space $\Bcal^{n+\ell}$ is $\dim(\Bcal^{n+\ell}) = n^2 +\ell$. Define 
\begin{alignat*}{8}
\Ab_k(i):=\begin{cases}
 A_k, & \text{for all } k \in \Ecal \cup \Ical, \ i=1 \\
 1,     & \text{if } k \in \Ical \text{ and } i=k-m+1.  \\
\end{cases}		
\end{alignat*}
Similarly, define 
\begin{alignat*}{8}
\Cb(i):=\begin{cases}
 C, &   i=1, \\
 0,    & \text{otherwise}.
\end{cases}		\qquad \qquad \qquad \qquad \ 
\end{alignat*}

Problem (\ref{sdp_primal1}) can be reformulated as a block diagonal semidefinite program in \textit{standard form}  as follows
\begin{alignat*} {6}
& \underset{\Xb \in \Bcal^{n+\ell}}{\text{minimize}}  
	 & & \Cb \bullet \Xb     \\
	 & \text{subject to} \ \
	&& \Ab_k \bullet \Xb  = b_k,  \ \  && \text{for all } k=1,\dots, m+\ell, \\
	&&&  \Xb \succeq 0,  
 \end{alignat*} 
where $\Cb:=\diag[C(1),\dots,C(\ell+1)] \in \Bcal^{n + \ell}$, $\Ab_k:=\diag[$ $\Ab_k(1),\dots,\Ab_k(\ell+1)] \in \Bcal^{n + \ell}$, $\Xb:= \diag[X(1),\dots,X(\ell+1)] \in \Bcal^{n+\ell}$. By $\Xb \succeq 0$ we mean $\Xb(i) \succeq 0$ for all $i=1,\dots,\ell+1.$ It is straightforward to see that for  $k \in \Ical$ and $i=k-m+1$, $\Xb(i) = 0$ if and only if $k \in \Acal(X). $

For each $i=1,\dots,\ell+1$, let $\Xb(i):=\Qb(i)\Lambdab(i) \Qb(i)^*$ be the singular value decomposition of $\Xb(i)$ and partition $\Qb(i):=\bmat{\Qb_1(i) & \Qb_2(i)},$ with $\Qb_1(i)$ and $\Qb_2(i)$ corresponding to the nonzero  and zero eigenvalues of $\Xb(i),$ respectively. When $\Xb(i) \in \R_+$ one of $\Qb_1(i)$, $\Qb_2(i)$ is the scalar 1, and the other is the empty matrix.  A direct application of Definition 3 in \cite{nayakkankuppam1999conditioning} leads to the definition of primal and dual nondegeneracy in (c.f. Definition \ref{def:nondegeneracy_primal}-\ref{def:nondegeneracy_dual}) in Section \ref{sec:IIC}.